\documentclass[12pt]{article}
\usepackage{amsfonts}
\usepackage{amssymb,amsmath,amstext}
\usepackage[mathscr]{eucal}

\begin{document}

\renewcommand{\labelenumi}{\theenumi)}
\newcommand{\qed}{\mbox{\raisebox{0.7ex}{\fbox{}}}}
\newtheorem{theorem}{Theorem}
\newtheorem{proposition}{Proposition }
\newtheorem{problem}{Problem}
\newtheorem{defin}{Definition}
\newtheorem{lemma}{Lemma}
\newtheorem{corollary}{Corollary}
\newtheorem{example}{\sc Example}
\newtheorem{remark}{\sc Remark}
\newtheorem{nt}{\sc Note}
\renewcommand{\thent}{}
\renewcommand{\theremark}{}
\newenvironment{pf}{\medskip\noindent{Proof:}
                \enspace}{\hfill \qed \newline \medskip}
\newenvironment{note}{\begin{nt} \em}{\end{nt}}
\newenvironment{exa}{\begin{example} \em}{\end{example}}
\newenvironment{defn}{\begin{defin} \em}{\end{defin}}
\newenvironment{rmk}{\begin{remark}
      \em}{\end{remark}}
\setlength{\unitlength}{12pt}
\newcommand{\comb}[2]{\mbox{$\left(\!\!\begin{array}{c}
            {#1} \\[-0.5ex] {#2} \end{array}\!\!\right)$}}

\title{Counting conjugacy classes of subgroups in a finitely generated
  group
\thanks{Supported by the RFBR (grant 03-01-00104) and by
INTAS (grant 03-51-3663) }}
\author
{Alexander Mednykh\\[1ex]
{\small\em Sobolev Institute of Mathematics, Novosibirsk State
University,}\\{\small\em  630090, Novosibirsk, Russia}\\{
\small\em    mednykh@math.nsc.ru }}

\date{}
\maketitle

\begin{abstract} A new general formula for  the number of conjugacy classes of
subgroups of given index in a finitely generated
  group is obtained.


\end{abstract}

\section{Introduction}
Let $M_{\Gamma}(n)$ denote the number of subgroups of index $n$
in a group $\Gamma,$ and $N_{\Gamma}(n)$ be the number of
conjugacy classes of such groups. The last function counts the
isomorphism classes of transitive permutation representations of
degree $n$ of $\Gamma$ and hence, also the equivalence classes of
$n$-fold unbranched connected coverings of a topological space
with fundamental group $\Gamma.$

 Each subgroup $K$ of index $n$ in a group $\Gamma$
determines a transitive action of degree $n$ of $\Gamma$ on the
cosets of $K,$ and each transitive action of $\Gamma$ on a set of
cardinality $n$ is isomorphic to such an action, with $K$
uniquely determined up to conjugacy in $\Gamma.$ Starting with
such a consideration M. Hall \cite{Hal49} determined the numbers
of subgroups $M_{\Gamma}(n)$ for a free group $\Gamma={\rm  F}_r$
of rank $r.$ Later V. Liskovets \cite{L71} developed  a new
method for calculation of $N_{\Gamma}(n)$ for the same group.
Both functions $M_{\Gamma}(n)$ and $N_{\Gamma}(n)$ for the
fundamental group $\Gamma$ of a closed surface were obtained in
\cite{M82} and \cite{MP86} for orientable and non-orientable
surfaces, respectively. See also \cite{M88} and \cite{JL} for the
case of the fundamental group  of the Klein bottle and a survey
\cite{KL01} for related problems. In all these cases the problem
of calculation of $M_{\Gamma}(n)$ was solved essentially due to
the ideas by Hurwitz and Frobenius to contribute the
representation theory of symmetric groups as the main tool
(\cite{H1},\,\cite{H2}). The solution for the problem to finding
$N_{\Gamma}(n)$ was based on the further  development of the
Liskovets method (\cite{L71},\,\cite{L98}). In \cite{LM1} and
\cite{LM2}, these ideas were applied to determine $M_{\Gamma}(n)$
for the fundamental groups of some Seifert spaces. Asymptotic
formulas for $M_{\Gamma}(n)$ in many important cases were
obtained in series of papers by T.~W.~ M\"uller and his
collaborators (\cite{Mue1},\cite{MueSh},\cite{MueP}). An
excellent exposition of the above results and related subjects is
given in the book \cite{LS}.

In the present paper, a new formula for the number of conjugacy
classes of subgroups of given index in a finitely generated group
is obtained. For the sake of simplicity, we require the group
$\Gamma$ to be finitely generated. But, indeed, all results of
the paper remain to be true for the groups having only finitely
many subgroups of each finite index (which is always  the case
for finitely generated groups).

The main counting principle suggested in  Section 2 of the paper
is rather universal. It can be applied to Fuchsian groups to
calculate the number of non-equivalent surface coverings (Section
3) as well as the number of unrooted maps on the surface
\cite{MN}. In Section 4, some general approach is developed to
find the number of non-equivalent unbranched coverings of a
manifold with finitely generated fundamental group.

\section{The main counting  principle}

 Denote by ${\rm Epi}(K,\,\texttt{Z}_\ell)$ the set of
epimorphisms of a group $K$ onto the cyclic group
$\texttt{Z}_\ell$ of order $\ell$ and by $|E|$ the cardinality of
a set $E.$

The main result of this paper is the following counting principle.

\begin{theorem} Let $\Gamma $ be a finitely generated group.
Then the number of  conjugacy classes of subgroups of index $n$
in the group $\Gamma $ is
  given by the formula
  $$N_{\Gamma}(n)=\frac{1}{n}\sum_{\substack{
   \ell|\,n \\
    \ell\,m=n } }\sum_{K\underset{m}{<}\,\Gamma }|{\rm Epi}(K,\,\texttt{Z}_\ell)|,$$
  where  the sum $\sum
  \limits_{K\underset{m}{<}\,\Gamma } $ is taken over all subgroups
  $K $ of index m in the group $\Gamma.$
\end{theorem}
\begin{pf} Let $P$ be a subgroup in $\Gamma$ and
$N(P,\Gamma)$ is the normalizer of $P$ in the group $\Gamma.$ We
need the following two elementary lemmas.

\begin{lemma}\label{number} The number of conjugacy classes of subgroups of index $n$ in the
group $\Gamma$ is given by the formula
$$
N_{\Gamma}(n)=\frac1n\sum\limits_{P\underset{n}{<}\,\Gamma}|N(P,\Gamma)/P|.
$$
\end{lemma}
\begin{pf}
Let $E$ be a conjugacy class of subgroups of  index $n$ in the
group $\Gamma.$ We claim   that

$$
\sum\limits_{P\in E}|N(P,\Gamma)/P|=n.
$$
Indeed, let $P'\in E.$ Then $|E|=|\Gamma :N(P',\Gamma)|$ and for
any $P\in E$ the groups $N(P,\Gamma)/P$ and $N(P',\Gamma)/P'$ are
isomorphic. We have
$$\sum\limits_{P\in E}|N(P,\Gamma)/P|=|E||N(P',\Gamma)/P'|=|\Gamma:N(P',\Gamma)||N(P',\Gamma):P'|=
|\Gamma:P'|=n.$$ Hence,
$$n\,N(n)=\sum\limits_E n=\sum\limits_E\sum\limits_{P\in E}|N(P,\Gamma)/P|=
\sum\limits_{P\underset{n}{<}\,\Gamma }|N(P,\Gamma)/P|,$$ where
the sum $\sum\limits_E$ is taken over all conjugacy classes $E$ of
subgroups of index $n$ in the group $\Gamma.$
\end{pf}
\bigskip
\begin{lemma}\label{normalizer}Let $P$ be a subgroup of index $n$ in the group
$\Gamma.$ Then
$$
|N(P,\Gamma)/P|=\sum\limits_{\substack{
   \ell|\,n \\
    \ell\,m=n } }
\sum\limits_{P\underset{Z_\ell}{\triangleleft}K\underset{m}{<}\,\Gamma
}\phi(\ell),$$ where $\phi(\ell)$ is the Euler
 function and the second sum is taken over all subgroups $K$ of
 index $m$ in $\Gamma$ containing $P$ as a normal subgroup with
 $K/P\cong\texttt{Z}_\ell.$ The sum vanishes if there are no such
 subgroups.
 \end{lemma}
 \begin{pf}Set $G=N(P,\Gamma)/P.$ Since $P\lhd N(P,\Gamma)<\Gamma$ and $P\underset{n}{<}\Gamma,$ the order
 of any cyclic subgroup of $G$ is a divisor of $n.$

 Note that there is a one-to-one correspondence between cyclic subgroups
 $\texttt{Z}_\ell$ in the group $G$ and subgroups $K$ satisfying
 $ P\underset{Z_\ell}{\triangleleft}K\underset{m}{<}\,\Gamma,$ where $\ell\,m=n.$

Given a cyclic subgroup $\texttt{Z}_\ell<G$ there  are exactly
$\phi(\ell)$ elements of $G$ which generate $\texttt{Z}_\ell.$

Hence,
$$
|G|=\sum\limits_{\ell|\,n}\phi(\ell) \sum\limits_{\texttt{Z}_\ell
< \,G}1=\sum\limits_{\ell|\,n}
\phi(\ell)\sum\limits_{P\underset{Z_\ell}{\triangleleft}K\underset{m}{<}\,\Gamma}1=\sum\limits_{\ell|\,n}
\sum\limits_{P\underset{Z_\ell}{\triangleleft}K\underset{m}{<}\,\Gamma
}\phi(\ell).
$$
 \end{pf}
 We finish the proof of the theorem by applying Lemma~\ref{number} and Lemma~\ref{normalizer} for $\ell\,m=n:$

$$\begin{array}{c}
  n\,N(n)=\sum\limits_{ P\underset{n}{<}\,\Gamma}|N(P,\Gamma)/P|= \\{}\\ \sum\limits_{P\underset{n}{<}\,\Gamma}\sum\limits_{\ell|\,n}
   \sum\limits_{P\underset{\,Z_\ell}{\triangleleft}\,K \underset{m}{<}\,\Gamma}\phi(\ell)=
   \sum\limits_{\ell|\,n}\sum\limits_{P\underset{n}{<}\,\Gamma}\sum\limits_{P\underset{Z_\ell}{\triangleleft}K\underset{m}{<}\,\Gamma}\phi(\ell)=\\{}\\
  \sum\limits_{\ell|\,n}\sum\limits_{K\underset{m}{<}\,\Gamma}\sum\limits_{P\underset{\,Z_\ell}{\triangleleft}\,K}\phi(\ell)
  =\sum\limits_{\ell|\,n}\sum\limits_{K\underset{m}{<}\,\Gamma}|{\rm Epi}(K,\texttt{Z}_\ell)|.
\end{array}$$
The last equality is a consequence of the following observation.
Given subgroup $P,~\,P~\underset{\,Z_\ell}{\triangleleft}~K $
there are exactly $\phi(\ell)$ homomorphisms
$\psi:K\to\texttt{Z}_\ell,$ with $\rm{Ker}\,(\psi)=P.$
\end{pf}

\bigskip
Denote by ${\rm Hom}(\Gamma,\,\texttt{Z}_\ell)$ the set of
homomorphisms of a group $\Gamma$ into the cyclic group
$\texttt{Z}_\ell$ of order $\ell.$ Since $|{\rm
Hom}(\Gamma,\,\texttt{Z}_\ell)|=\sum\limits_{d|\,\ell}|{\rm
Epi}(\Gamma,\,\texttt{Z}_d)|,$ by the M\"obius inversion formula
we have the following result

\bigskip
\begin{lemma}\label{Gareth}({\bf G.\,Jones}\,\cite{J95})
$$|{\rm Epi}(\Gamma,\,\texttt{Z}_\ell)|=\sum\limits_{d|\,\ell}\mu(\frac{\ell}{d})|{\rm Hom}(\Gamma,\,\texttt{Z}_d)|,$$
where $\mu(n)$ is the M\"obius function.
\end{lemma} This lemma essentially
simplifies the calculation of $|{\rm
Epi}(\Gamma,\,\texttt{Z}_\ell)|$ for a finitely generated group
$\Gamma.$ Indeed, let ${\rm
H}_1(\Gamma)=\Gamma/[\Gamma,\,\Gamma]$ be the first homology
group  of  $\Gamma.$ Suppose that $\rm{
H}_1(\Gamma)=\texttt{Z}_{m_1}\oplus\texttt{Z}_{m_2}\oplus\ldots\oplus\texttt{Z}_{m_s}\oplus\texttt{Z}^r.$
Then we have
\bigskip
\begin{lemma}\label{Derevnin}
$$|{\rm Epi}(\Gamma,\,\texttt{Z}_\ell)|=\sum\limits_{d|\,\ell}\mu(\frac{\ell}{d})\,(m_1,d)\,(m_2,d)\ldots(m_s,d)\,d^r,$$
where $(m,d)$ is the greatest common divisor of $m$ and $d.$
\end{lemma}
\begin{pf}Note that $|{\rm Hom}(
\texttt{Z}_m,\,\texttt{Z}_d)|=(m,d)$ and $|{\rm
Hom}(\texttt{Z},\,\texttt{Z}_d)|=d.$ Since the group
$\texttt{Z}_d$ is Abelian, we obtain
$$
\begin{array}{c}
 |{\rm Hom}(\Gamma,\,\texttt{Z}_d)|= |{\rm Hom}({\rm H}_1(\Gamma),\,\texttt{Z}_d)|= (m_1,d)\,(m_2,d)\ldots(m_s,d)\,d^r.\\
 {}
\end{array}
$$ Then the result follows from Lemma~\ref{Gareth}.
\end{pf}

In particular, we have
\bigskip

\begin{corollary}\label{surface}
\begin{eqnarray*}
\text{(i)} & &  \text{Let}\ \rm{F}_r \, \text{be a free group of
rank} \  r. \ \text{Then} \
\rm{H}_1({\rm  F}_r)=\texttt{Z}^{\,r} \,\,\text{and}\  \\
 & & |{\rm Epi}({\rm  F}_r,\,\texttt{Z}_\ell)|~=~\sum\limits_{d|\,\ell}\mu(\frac{\ell}{d})\,d^r. \nonumber \\
\text{(ii)} & & \,\text{Let}\ \Phi_g=<a_1,b_1,\ldots,a_g,b_g:
\prod\limits_{i=1}^{g}[a_i,\,b_i]=1> \text{be the fundamental} \\
  & &  \text{group of a closed orientable surface of genus} \ g. \
  \text{Then}\
\rm{H}_1(\Phi_g)=\texttt{Z}^{\,{2g}} \\
  & & \text{and}\ \,\, |{\rm Epi}(\Phi_r,\,\texttt{Z}_\ell)|=\sum\limits_{d|\,\ell}\mu(\frac{\ell}{d})\,d^{2g}.\\
\text{(iii)}  & & \ \text{Let} \ \Lambda_p=<a_1,a_2,\ldots,a_p:
\prod\limits_{i=1}^{p}a_i^2=1> \ \text{be the fundamental group} \\
  & & \text{of a closed non-orientable surface of genus} \  p. \
  \text{Then}\
\rm{H}_1(\Lambda_p)=\texttt{Z}_{2}\oplus\texttt{Z}^{\,{p-1}} \\
  & & \text{and}\
 \,\,|{\rm Epi}(\Lambda_p,\,\texttt{Z}_\ell)|=\sum\limits_{d|\,\ell}\mu(\frac{\ell}{d})\,(2,d)\,d^{p-1}.
\end{eqnarray*}

\end{corollary}
\bigskip

Note that the fundamental group of any compact surface (orientable
or not, possibly, with non-empty boundary) is one of the three
groups ${\rm  F}_r,\,\Phi_g$  and $\Lambda_p$ listed in Corollary
~\ref{surface}. In the next two sections we identify
 the number of conjugacy
classes of  subgroups of index $n$  in the group $\Gamma$ and the
number of equivalence classes of $n$-fold unbranched connected
coverings of a  manifold with fundamental group $\Gamma.$
   \section{Counting surface coverings}

Recall that the fundamental group $\pi_1(\mathcal{B})$ of a
bordered surface $\mathcal{B}$ of Euler characteristic
$\chi=1-r,\,r\ge 0,$ is a free group ${\rm  F}_r$ of rank $r.$ An
example of such a surface is the disc $\mathcal{D}_r $ with $r$
holes.  As the first application of the counting principle
(Theorem 1) we have the following result obtained earlier by
V.\,Liskovets \cite{L71}
\bigskip
\begin{theorem}\label{Liskovets} Let $\mathcal{B}$ be a bordered surface with the
fundamental group $\pi_1(\mathcal{B})={\rm  F}_r.$ Then the
number of non-equivalent $n-$fold coverings of $\mathcal{B}$ is
  given by the formula
  $$N(n)=\frac{1}{n}\sum\limits_{\substack{
   \ell|\,n \\
    \ell\,m=n
    }}\sum\limits_{d|\,\ell}\mu(\frac{\ell}{d})\,d^{(r-1)m+1}M(m),
  $$where $M(m)$ is the number of subgroups of index
  $m$ in the group ${\rm  F}_r .$
\end{theorem}
\begin{pf} Note that all subgroups of index $m$ in ${\rm  F}_r$ are
isomorphic to $\Gamma_m={\rm  F}_{(r-1)m+1}.$ By Theorem 1 we have
$$N(n)=\frac{1}{n}\sum\limits_{\substack{
   \ell|\,n \\
    \ell\,m=n
    }}|{\rm Epi}(\Gamma_m,\, \texttt{Z}_\ell)|\cdot M(m),$$

    By Corollary 1(i) we get
$$|{\rm Epi}(\Gamma_m,\, \texttt{Z}_\ell)|=\sum\limits_{d|\,\ell}\mu(\frac{\ell}{d})\,d^{(r-1)m+1}$$
and the result follows.
\end{pf}
By the M.\,Hall  recursive formula \cite{Hal49} the number of
subgroups of index
  $m$ in the group ${\rm  F}_r $ is equal to
  $M(m)=\displaystyle{\frac{t_{m,\,r}}{(m-1)!}},$ where
  $$t_{m,\,r}={m!}^r-\sum_{j=1}^{m-1}\comb{m-1}{j-1}{(m-j)!}^r\,t_{j,\,r},\,\,\,t_{1,\,r}=1.$$

\bigskip
 The next result was obtained in \cite{M82} in a rather complicated way
\begin{theorem}\label{Mednykh} Let $\mathcal{S}$ be a  closed orientable  surface with the
fundamental group $\pi_1(\mathcal{S})=\Phi_g.$ Then the number of
non-equivalent $n-$fold coverings of $\mathcal{S}$ is
  given by the formula
  $$N(n)=\frac{1}{n}\sum\limits_{\substack{
   \ell|\,n \\
    \ell\,m=n
    }}\sum\limits_{d|\,\ell}\mu(\frac{\ell}{d})\,d^{2(g-1)m+2}M(m),
  $$where $M(m)$ is the number of subgroups of index
  $m$ in the group $\Phi_g .$
\end{theorem}
\begin{pf} In this case, by the Riemann-Hurwitz formula all subgroups of index $m$ in $\Phi_g$ are
isomorphic to the group $K_m=\Phi_{(g-1)m+1}.$ By  the main
counting principle we have
$$N(n)=\frac{1}{n}\sum\limits_{\substack{
   \ell|\,n \\
    \ell\,m=n
    }}|{\rm Epi}(K_m,\, \texttt{Z}_\ell)|\cdot M(m),$$where
$$|{\rm Epi}(K_m,\, \texttt{Z}_\ell)|=\sum\limits_{d|\,\ell}\mu(\frac{\ell}{d})\,d^{2(g-1)m+2}$$
is given by Corollary 1(ii) .
\end{pf}

Let $\mathcal{N}$ be a closed non-orientable surface of genus $p$
with the fundamental group $\pi_1(\mathcal{N})=\Lambda_p.$ Denote
by $\mathcal{N}_m^+$ and $\mathcal{N}_m^-$ an orientable and
non-orientable $m-$fold coverings of $\mathcal{N},$ respectively
and set $\Gamma_m^{+}=\pi_1(\mathcal{N}_{m}^{+})$ and
$\Gamma_m^{-}=\pi_1(\mathcal{N}_{m}^{-}).$  For simplicity, we
will refer to $\Gamma_m^{+}$ and $\Gamma_m^{-}$ as orientable and
non-orientable subgroups of index $m$ in $\Lambda_p,$
respectively. By
  the Riemann-Hurwitz formula we  get
  $$2\,{\rm genus\,}(\mathcal{N}_m^{+})-2=m(p-2)\,\,\,{\rm and}\,\,\,{\rm
  genus\,}(\mathcal{N}_m^{-})-2=m(p-2),$$ where $p={\rm genus\,}(\mathcal{N}).$
Hence $\Gamma_m^{+}=\Phi_{\frac m2(p-2)+1}$ and $\Gamma_m^{-}=
\Lambda_{m(p-2)+2}.$

By the main counting principle, the number of non-equivalent
$n-$fold coverings of $N$ is
  given by the formula
  $$N(n)=\frac{1}{n}\sum\limits_{\substack{
   \ell|\,n \\
    \ell\,m=n
    }}(|{\rm Epi}(\Gamma_m^{+},\, \texttt{Z}_\ell)| \cdot M^{+}(m)
  +|{\rm Epi}(\Gamma_m^{-},\, \texttt{Z}_\ell)| \cdot M^{-}(m)),$$where
  $M^{+}(m)$ and $M^{-}(m)$ are the numbers of orientable and non-orientable subgroups of index
  $m$ in the group $\Lambda_p ,$ respectively.

   By Corollary 1(ii) and Corollary 1(iii), we have
$$|{\rm Epi}(\Gamma_m^{+},\,
\texttt{Z}_\ell)|=\sum\limits_{d|\,\ell}\mu(\frac{\ell}{d})\,d^{m(p-2)+2
}\,\,\,{\rm and}\,\,\,
  |{\rm Epi}(\Gamma_m^{-},\, \texttt{Z}_\ell)|=\sum\limits_{d|\,\ell}\mu(\frac{\ell}{d})\,(2,d)\,d^{m(p-2)+1
  }.$$

  As a result, we have proved the following theorem obtained
  earlier in \cite{MP86} by making use of a cumbersome combinatorial technique.
  \bigskip
  \begin{theorem}\label{MP} Let $n$ be a  closed orientable surface with the
fundamental group $\pi_1(\mathcal{N})=\Lambda_p.$ Then the number
of non-equivalent $n-$fold coverings of $\mathcal{N}$ is
  given by the formula
  $$N(n)=\frac{1}{n}\sum\limits_{\substack{
   \ell|\,n \\
    \ell\,m=n
    }}\sum\limits_{d|\,\ell}\mu(\frac{\ell}{d})\,(d^{m(p-2)+2}M^{+}(m)+(2,d)\,d^{m(p-2)+1}M^{-}(m)),
  $$where
  $M^{+}(m)$ and $M^{-}(m)$ are the numbers of orientable and non-orientable subgroups of index
  $m$ in the group $\Lambda_p ,$ respectively.
\end{theorem}
For completeness note that (\cite{M82},\,\cite{MP86})  if
$\Gamma=\Phi_g$ or $\Lambda_p$ then the number $M(m)$ of
subgroups of index
  $m$ in the group $\Gamma $ is equal to
  $$R_\nu(m)=m\sum_{s=1}^m\frac{(-1)^{s+1}}{s}\sum_{\substack{
    i_1+i_2+\ldots+i_s=m\\
    i_1,\,i_2,\,\ldots\,,i_s\ge 1 }
    }\beta_{i_1}\beta_{i_2}\cdots\beta_{i_s},$$ where
    $\beta_k=\sum\limits_{\chi\in D_k}\displaystyle{\big(\frac{k!}{f^{\chi}}\big)^\nu},\,\,\,D_k $ is the set of
    irreducible representations of a symmetric group $S_k,
    \,\,f^{\chi}$ is the degree of the representation $\chi,$
     $\nu=2g-2$ for $\Gamma=\Phi_g$ and  $\nu=p-2$ for $\Gamma=\Lambda_p.$ Moreover, in the
     latter case,
     $M^{+}(m)=0$ if $m$ is odd, $M^{+}(m)=R_{2\,\nu}(\frac m2)$ if $m$ is even, and
     $M^{-}(m)=M(m)-M^{+}(m).$ Also, the number of
     subgroups can be found by the following recursive formula
     $$M(m)=m\,\beta_m-\sum_{j=1}^{m-1}\beta_{m-j}M(j),\,\,\,M(1)=1.$$

 \section{Non-equivalent coverings of manifolds }

All manifolds in this section are supposed to be connected, with
finitely generated fundamental group. No restriction on dimension
is given. The manifolds under consideration can be closed, open or
bordered, orientable or not. The  following theorem is just a
topological version of Theorem 1.
\begin{theorem}\label{homotopy}   Let $\texttt{M}$ be a connected manifold with finitely generated
fundamental group $\Gamma=\pi_1(\texttt{M}).$  Then the number of
non-equivalent $n-$fold coverings of $\texttt{M}$ is
  given by the formula
  $$N(n)=\frac{1}{n}\sum_{\substack{
   \ell\,|\,n \\
    \ell\,m=n } }\sum\limits_{\displaystyle{\Phi \in \,\mathfrak{F}_m }}
  |{\rm Epi}( \Phi ,\, \texttt{Z}_\ell)|\cdot M_{\Phi,\,\Gamma}(m),$$
where $\mathfrak{F}_m $ is the set of groups arising as
fundamental groups of
  $m-$fold coverings  of $\texttt{M}$ and $M_{\Phi,\,\Gamma}(m) $
  is the number of  subgroups of index
  $m$ in the group  $\Gamma$ that are isomorphic to $ \Phi.$
\end{theorem}
Taking into
  account  Lemma~\ref{Gareth}, we obtain the following result as
  an immediate corollary of Theorem~\ref{homotopy}
  \bigskip
\begin{theorem}\label{homotophom}   Let $\texttt{M}$ be a connected manifold with finitely generated
fundamental group $\Gamma=\pi_1(\texttt{M}).$  Then the number of
non-equivalent $n-$fold coverings of $\texttt{M}$ is
  given by the formula
  $$N(n)=\frac{1}{n}\sum_{\substack{
   \ell\,|\,n \\
    \ell\,m=n } }\sum\limits_{\displaystyle{\Phi \in \,\mathfrak{F}_m }}
 \sum\limits_{d|\ell}\mu(\frac{\ell}{d})
  |{\rm Hom}({\Phi} ,\, \texttt{Z}_d)|
  \cdot M_{\Phi,\,\Gamma}(m),$$
where $\mathfrak{F}_m $ is the set of groups arising as
fundamental groups of
  $m-fold$ coverings  of $\texttt{M}$ and $M_{\Phi,\,\Gamma}(m) $
  is the number of  subgroups of index
  $m$ in the group  $\Gamma$ that are isomorphic to $ \Phi.$
\end{theorem}

\bigskip Let $\texttt{M}$ be a manifold and $\Gamma=\pi_1(\texttt{M}).$
Denote by ${\rm H}_1(\texttt{M})$ and  ${\rm H}_1(\Gamma)$ the
first homology over $\texttt{Z}$ of the manifold $\texttt{M}$ and
 the group $\Gamma,$  respectively.  Recall that ${\rm
H}_1(\Gamma)=\Gamma/[\Gamma,\,\Gamma]$ and ${\rm
H}_1(\texttt{M})={\rm H}_1(\Gamma).$

  Since the group  $ \texttt{Z}_d$ is Abelian, there is a one-to-one correspondence between the sets
  ${\rm Hom}(\Gamma ,\, \texttt{Z}_d)$
  and ${\rm Hom}(\rm {H}_1(\Gamma ),\, \texttt{Z}_d).$ As a result we
  have the following homological version of the previous theorem
  \bigskip
  \begin{theorem}\label{homology}   Let $\texttt{M}$ be a connected manifold with finitely generated
fundamental group $\Gamma=\pi_1(\texttt{M}).$  Then the number of
non-equivalent $n-$fold coverings of $\texttt{M}$ is
  given by the formula

  $$N(n)=\frac{1}{n}\sum_{\substack{
   \ell\,|\,n \\
    \ell\,m=n } }\sum\limits_{\displaystyle{\rm{H} \in \mathfrak{H}_m }}
  \sum\limits_{d|\ell}\mu(\frac{\ell}{d})
  |{\rm Hom}({\rm H} ,\, \texttt{Z}_d)|
  \cdot M^{'}_{\rm{H},\,\Gamma}(m),$$
where $\mathfrak{H}_m $ is the set of groups arising as
homologies
 of $m-$fold coverings  of $\texttt{M}$ and $M'_{\rm{H},\,\Gamma}(m) $
  is the number of  subgroups $F$ of index
  $m$ in the group  $\Gamma,$ with ${\rm H}_1(F)$ isomorphic to $ \rm{H}.$

\end{theorem}
 {\bf Acknowledgement:}
 % The authors are very appreciated
%to Prof. V.~A.~Liskovets for helpful remarks and suggestions. The
The author gratefully thanks Com$^2$MaC, Pohang University of
Science and Technology (Korea) for hospitality during preparation
of the present paper.

\end{document}